\documentclass[12pt]{amsart}
\usepackage{graphicx,color,amssymb,latexsym,amsfonts,amsmath}
\usepackage{graphicx}
\usepackage{graphics}

\newtheorem{thm}{Theorem}[section]

\theoremstyle{definition}

\theoremstyle{remark}

\numberwithin{equation}{section}

\begin{document}

\title[Prediction of missing observations]{Prediction of missing observations
  by a control method }%
\author{Vyacheslav M. Abramov}%
\author{Fima C. Klebaner}
\address{School of Mathematical Sciences, Monash University, Building
28M, Clayton campus, Clayton, VIC 3800, Australia}%
\email{vyacheslav.abramov@sci.monash.edu.au}%
\email{fima.klebaner@sci.monash.edu.au}%


\begin{abstract} Consider  a time series with missing observations
but a known   final point. Using control theory ideas we
estimate/predict these missing observations.  We obtain recurrence
equations which minimize sum of squares of a control sequence. An
advantage of this method is in easily computable formulae and
flexibility of its application to different structures of missing
data.
\end{abstract}
\maketitle
\section{Introduction}\label{Introduction}

Analysis and forecasting missing data is a well-known area of
statistics going back to   earlier works of Bartlett
\cite{Bartlett (1)}, Tocher \cite{Tocher (1)}, Wilks \cite{Wilks
(1)}, Yates \cite{Yates (1)} and many others (see review paper
\cite{Afifi and Elashoff (1)}). There is a large number of review
papers and books related to this subject, \cite{Afifi and Elashoff
(1)}, \cite{Alison (1)}, \cite{Jones (1)}, \cite{Little (1)} and
\cite{Little and Rubin (1)}, to mention a few. There are various
approaches to missing data,  including Bayes methods \cite{Carlin
Louis (1)}, maximum likelihood, multiple imputations methods,
methods of non-parametric regression and others, e.g. \cite{Alison
(1)}, \cite{David Little Samuhel and Triest (1)}, \cite{Wang
Sedransk and Jinn (1)}.

In the present paper we suggest a new method of predicting a
special class of missing observations in different time series
including regression and auto-regression models. We suggest a
simple recurrence procedure, and to the authors knowledge, it is
new and simpler than the computational procedures that were known
before.

We study autoregressive time series with missing observations,
which we propose to predict using a control method. This method is
developed for different types of autoregressive models including
AR($p$) models in the case of scalar variables and AR(1) in the
case of vector-valued observations. Forecasting missing data in
autoregressive time series has received a special attention in the
literature: \cite{Jones (1)}, \cite{Harvey and Pierse (1)},
 \cite{Kharin and Huryn (2)},
\cite{Kohn and Ansley (1)}, \cite{Little and Rubin (1)},
\cite{Pourhamadi (1)}, \cite{Sargan and Drettakis (1)} and
\cite{Stoffer (1)}. The typical approach for forecasting missing
data in
 autoregressive models considered in  most of papers is based
 on maximization of likelihood ratio, which can be computationally intensive.
 The approach of the present
 paper, referred to as a \textit{control method} allows to obtain   easily
 computable formulae for missing data.


It is known that a one-dimensional AR($p$) model can be
transformed to the special case of a $p$-dimensional AR(1) model
(e.g. Anderson \cite{Anderson (1)}). In the present paper we
consider both one-dimensional AR($p$) models and multidimensional
AR(1) models nevertheless. The representations obtained in the
case of one-dimensional AR($p$) models are simpler for
computations   than that for a multidimensional AR(1) model.
Whereas   representations for a one-dimensional AR($p$) model is
  recurrence formulae and can be calculated directly, the
computations for a multidimensional AR(1) model requires two
steps. In the first step we calculate the vector norm, and then
the vector corresponding to a missing value.

We assume that in the   time series:
\begin{equation}\label{I1}
\begin{array}{llllllllll}
&x_1, &x_2, &\ldots, &x_{n_0}, &\widetilde{x}_{n_0+1},
&\widetilde{x}_{n_0+2}, &\ldots, &\widetilde{x}_{N-1}, &
x_{N}=\overline{x}
\end{array}
\end{equation}
  the first $n_0$ observations
  are   known/ observed as well as the
last value $x_N=\overline{x}$ is assumed to be given too. The
values $\widetilde{x}_{n_0+1}$, $\widetilde{x}_{n_0+2}$,\ldots,
$\widetilde{x}_{N-1}$  are missing.

This set up   may have various applications, for example in
economics and finance, where historical data indices are given,
while the last value can be obtained from financial derivatives,
or might be set externally. In finance, for example, on basis of
historical volatilities and a future value obtained from options
one predicts the dynamics of the volatility.

Although the paper concerns with data structure \eqref{I1}, the
results   can be extended to different more complicated structures
of missing values. Indeed, consider for instance the following
data
\begin{equation}\label{I2}
\begin{array}{llllllll}
&x_1, &\ldots, &x_{n_0}, &\widetilde{x}_{n_0+1},  &\ldots,
&\widetilde{x}_{n_1},\\
&x_{n_1+1},&\ldots,&x_{n_2},&\widetilde
x_{n_2+1},&\ldots,&\widetilde x_{N-1}, &x_{N}=\overline{x}.
\end{array}
\end{equation}
Here in \eqref{I2} the data indexed from 1 to $n_0$ and from
$n_1+1$ to $n_2$ are known, the last point $x_N=\overline{x}$ is
assigned, and the rest of data are missing. Then we have two
groups of missing data, and standard decomposition arguments can
be used to reduce analysis to that of a single group with missing
data.

The paper is structured as follows. In Section \ref{Approaches} we
discuss the known methods of forecasting missing observation as
well as the method of the present paper with comparisons. In the
following section we discuss forecasting missing data by a control
method in order of increasing complexity. Specifically, in Section
\ref{autoregressive} we study the problem for the simplest AR(1)
model of time series, and in Section \ref{AR models} we extend the
results for AR($p$) models, $p\geq1$. The multi-dimensional
observations of AR(1) model are studied in Section
\ref{Multi-autoregressive}. Then, in Section
\ref{Multi-regression} the problem is solved for models of
regression. The results of this section are easily understandable
and simple.
In Section \ref{Numerical work} two numerical examples are
considered in finance and  archaeology.

\section{Review of methods for missing values}\label{Approaches}

There is a large number of papers on  estimation and forecasting
of missing observation  in autoregressive models.

Jones \cite{Jones (1)} provides the method for calculation of
exact likelihood function of stationary ARMA time series. The
method is based on Akake's Markovian representation and
application of Kalman's recursions \cite{Brockwell and Davis (1)}.
 An advantage of Kalman's recursions is that the matrices
and vectors being used in calculations have dimensions
$\max\{p,q+1\}$, where $p$ is the order of the auto-regression,
and $q$ is the order of moving average, rather than dimensions
corresponding to the number of observations.  
A non-linear optimization program is then used
 to find the maximum likelihood estimates.

Kohn and Ansley \cite{Kohn and Ansley (1)} study interpolation
missing data for non-stationary ARIMA models. The likelihood ratio
for these models does not exist in the  usual sense, and the
authors define marginal likelihood ratio. They show, that marginal
likelihood approach reduces in some cases to the usual likelihood
approach
Forecasting missing observations in   is based on a modified
Kalman filter, which has been introduced in the earlier paper of
these authors \cite{Kohn and Ansley (2)}.

Shin and Pantula \cite{Shin and Pantula (1)} discuss the testing
problem for a unit root in an autoregressive model where data are
available for each $m$-th period. The idea  is to use
characteristic polynomials and   properties of their coefficients.
Under special assumption \cite{Shin and Pantula (1)} estimate
parameters of ARMA($p$, $p$-1) by fitting an ARMA(1, $p$-1) model.
By using a Monte Carlo simulation, the results were compared by
those obtained in earlier papers of Pantula and Hall \cite{Pantula
and Hall (1)}, Said and Dickey \cite{Said and Dickey (1)} and Shin
and Fuller \cite{Shin and Fuller (1)}, who also studied the same
testing problem.

Forecasting in autoregressive models has also been studied by
Kharin and Huryn  \cite{Kharin and Huryn (2)} and \cite{Kharin and
Huryn (3)}.    \cite{Kharin and Huryn (2)}  investigate  the case
of unknown parameters of an autoregressive model based on the
so-called ``plug-in" approach. The ``plug-in" approach consists of
two steps: (i) estimation of the model parameters by some known
approach and (ii) forecasting,   based on estimation of the
parameters in the first step. This method has lower computational
complexity than other   methods, such as   straightforward joint
maximum likelihood estimation of the parameters and future values
of time series, or Expectation-Maximization algorithm (e.g. Little
and Rubin \cite{Little and Rubin (1)}, Jordan and Jacobs
\cite{Jordan and Jacobs (1)}). In \cite{Kharin and Huryn (3)} the
mean-squared error of maximum likelihood forecasting in the case
of missing values is obtained  for many autoregressive time
series.

The above-mentioned papers \cite{Kharin and Huryn (2)} and
\cite{Kharin and Huryn (3)} all study a general scheme of missing
data. Together with vector-valued time series they introduce a
binary vector characterizing a ``missing pattern" but the solution
for this general formulation is hard to implement  in practice.


 The aim of the present paper is
prediction (interpolation) of missing observations whereas the aim
of two above-mentioned papers is forecasting \textit{in the
presence} of missing observations, i.e. the forecasting procedure
takes into account  missing observations. Furthermore, the
approach of the present paper deals with   specific data
 structures (Section
\ref{Introduction}), and can be extended   to more complicated
structures of missing data. In the initial step we use   least
squares predictors for the preliminary extrapolation of missing
values. Then, taking into account the last known observation we
 make corrections by formulating and solving a
control problem. The   control problem is formulated in terms of
minimization of  sums of squares of errors, which in itself is
  a classical approach. However, our method of
is based on a novel application of the Cauchy-Schwartz inequality
in a simple case, and then   extended to other more complicated
cases. The use of the Cauchy-Schwartz inequality is a known
technique in optimization, e.g. \cite{Chigansky Liptser and
Bobrovsky (1)} and \cite{Iacus and Kutoyants (1)}, however, in the
context of prediction of missing data this method seems to be new.
In addition, this method yields easily  computable recurrence
formulae for missing values.

\section{A control method for missing data}\label{autoregressive}
In this section we consider autoregressive time series of the
following type:
\begin{equation}\label{ATS1}
\begin{array}{llllllllll}
&x_1, &x_2, &\ldots, &x_{n_0}, &\widetilde{x}_{n_0+1},
&\widetilde{x}_{n_0+2}, &\ldots, &\widetilde{x}_{N-1}, &
x_{N}=\overline{x}.
\end{array}
\end{equation}
The values $x_1$, $x_2$, \ldots, $x_{n_0}$ are assumed to be
observed, while by $\widetilde{x}_{n_0+1}$,
$\widetilde{x}_{n_0+2}$, \ldots, $\widetilde{x}_{N-1}$ we denote
estimates of missing observations.  The value ${x}_N=\overline{x}$
is also    known.
  It is convenient to denote this value   by
tilde, i.e. ${x}_N=\widetilde{x}_N=\overline{x}$.


\begin{thm}\label{thm1}
Best predictors for the missing values are given by
\begin{equation}\label{ATS10}
\widetilde x_n=a\widetilde x_{n-1}+b+\frac{ \overline{x}- \widehat
x_N}{\sum_{i=n_0+1}^Na^{2(N-i)}}\cdot a^{N-n},\;\;n=n_0+1,
  \ldots, N-1,
\end{equation}
where  the coefficients $a$ and $b$  are the least squares
solutions  of the autoregressive equations
\begin{equation*}
x_n=a x_{n-1}+b,
\end{equation*}
for  the first $n_0$ observations, $n=1,2,\ldots, n_0$.
\end{thm}

\begin{proof}
Taking into consideration   the first $n_0$ observed values one
can build the linear least square predictor as
\begin{equation}
\label{ATS2} \widehat x_n=a \widehat x_{n-1}+b
\end{equation}
for $n=n_0+1, n_0+2,\ldots, N$, where parameters $a$ and $b$ are
the regression coefficients. These $a$ and $b$ are then used for
control problem, which is to find the   unknown points, minimizing
sum of squares of controls leading to the known final value.
Namely,   for $n=n_0+1,\ldots,N$ ($\widetilde x_{n_0}$=$\widehat
x_{n_0}$)
\begin{equation}
\widetilde x_n=a \widetilde x_{n-1}+b+u_n.\label{ATS3}
\end{equation}
 It can be seen as a \textit{correction} of the initial
linear equation for $\widehat x_n$ with a control sequence $u_n$.
 The control problem is to minimize the sum of squares of
controls under the condition that the auto-regression ends up at
the specified point $\widetilde x_N$=$\overline{x}$
\begin{equation}
\min_{u_n:\;\widetilde x_N=\overline{x}} \ \sum_{n=n_0+1}^N
u_n^2.\label{ATS4}
\end{equation}
This minimization problem is solved as follows. By \eqref{ATS2}
and \eqref{ATS3}
\begin{equation}\label{ATS5}
u_n^2=(\widehat x_n-\widetilde x_n)^2,
\end{equation}
and taking into account that
\begin{equation*}
\widetilde x_N=\widetilde
x_{n_0}a^{N-n_0}+\sum_{n=n_0+1}^Na^{N-n}(b+u_{n})
\end{equation*}
and
\begin{equation*}
\widehat x_N=\widehat x_{n_0}a^{N-n_0}+b\sum_{n=n_0+1}^Na^{N-n},
\end{equation*}
from \eqref{ATS5} we obtain
\begin{equation}\label{ATS6}
u_N^{2}=\left(\sum_{n=n_0+1}^Na^{N-n}u_{n}\right)^{2}.
\end{equation}
By the Cauchy-Schwartz inequality,
\begin{equation}\label{ATS7}
\left(\sum_{n=n_0+1}^Na^{N-n}u_{n}\right)^{2}\leq
\sum_{n=n_0+1}^Na^{2(N-n)}\cdot\sum_{n=n_0+1}^Nu_{n}^{2}.
\end{equation}
The equality in \eqref{ATS7} is achieved if and only if
$a^{N-n}=cu_n$ for some constant $c$, and since the equality in
\eqref{ATS7} is associated with the minimum of the left-hand side
of \eqref{ATS7}, the problem reduces to find an appropriate value
$c=c^*$ such that
\begin{equation*}
u_n=c^*a^{N-n}.
\end{equation*}
Therefore,
\begin{equation*}
\widetilde x_N= \widehat x_N+c^*\sum_{i=n_0+1}^Na^{2(N-i)},
\end{equation*}
and then finally for $c^*$ we have:
\begin{equation}\label{ATS8}
c^*=\frac{\widetilde x_N- \widehat
x_N}{\sum_{i=n_0+1}^Na^{2(N-i)}}.
\end{equation}
Thus, the sequence $u_n$ satisfying \eqref{ATS4} is
\begin{equation*}
u_n=\frac{ \overline{x}- \widehat
x_N}{\sum_{i=n_0+1}^Na^{2(N-i)}}\cdot a^{N-n},
\end{equation*}
and its substitution for \eqref{ATS3} yields the desired result
\eqref{ATS10}.
\end{proof}

\section{Extension of the result for AR($p$) model}\label{AR
models} Under the assumption that \eqref{ATS1} is given, we first
find the best linear predictor for AR(2) model as
\begin{equation}
\label{AR2}
\widehat{x}_n=a_1\widehat{x}_{n-1}+a_2\widehat{x}_{n-2}+b,
\end{equation}
($n=n_0+2$, $n_0+3$+\ldots, $N$), and then extend the result for
the general case of AR($p$) model.

\begin{thm}\label{thm2}
For AR(2) model, the best predictor is given by
\begin{equation}\label{AR9}
\widetilde x_n=a_1\widetilde x_{n-1}+a_2\widetilde
x_{n-2}+b+\frac{ \overline{x}- \widehat
x_N}{\sum_{i=n_0+2}^N\gamma_{N-i}^2}\cdot \gamma_{N-n},
\end{equation}
where $n=n_0+2, n_0+3,\ldots, N-1$, and the coefficients
$\gamma_n$ are as follows: $\gamma_n=\alpha_n+\beta_n$,
\begin{eqnarray*}
\alpha_0 &=& 1,\nonumber\\
\beta_0 &=& 0,\nonumber\\
\alpha_n &=& a_1(\alpha_{n-1}+\beta_{n-1}) \ \ (n\geq1),\\
\beta_1 &=& 1,\nonumber\\
\beta_n &=& a_2(\alpha_{n-2}+\beta_{n-2})+1 \ \ (n\geq2)\nonumber.
\end{eqnarray*}

The coefficients $a_1$, $a_2$ and $b$ for equation \eqref{AR9} are
the minimum in the least-square sense of the autoregressive
equation
$$
x_n=a_1x_{n-1}+a_2x_{n-2}+b,
$$
which are obtained by the first $n_0$ observations.
\end{thm}

\begin{proof}
In the case of AR(2) model we have
\begin{equation}\label{AR3}
\widehat{x}_{n_0+2}=a_1\widehat{x}_{n_0+1}+a_2\widehat{x}_{n_0}+b,
\end{equation}
and similarly to \eqref{ATS3},
\begin{equation}
\label{AR4}
\widetilde{x}_n=a_1\widetilde{x}_{n-1}+a_2\widetilde{x}_{n-2}+b+u_n.
\end{equation}
($n=n_0+2, n_0+2,\ldots,N$).

Let us now consider the difference
$\widetilde{x}_N-\widehat{x}_N=u_N$. For this difference we have
the following expansion
\begin{equation}\label{AR5}
u_N=\sum_{n=n_0+2}^N \gamma_{N-n}u_n,
\end{equation}
with some coefficients $\gamma_{N-n}$. Now, the main task is to
determine these coefficients. Write $\gamma_n=\alpha_n+\beta_n$.
Then, using induction we obtain
\begin{eqnarray}\label{AR6}
\alpha_0 &=& 1,\nonumber\\
\beta_0 &=& 0,\nonumber\\
\alpha_n &=& a_1(\alpha_{n-1}+\beta_{n-1}) \ \ (n\geq1),\\
\beta_1 &=& 1,\nonumber\\
\beta_n &=& a_2(\alpha_{n-2}+\beta_{n-2})+1 \ \ (n\geq2)\nonumber.
\end{eqnarray}
Specifically, for the first steps we have the following. Setting
$N=n_0+2$ leads to the obvious identity
$u_N=\gamma_0u_N=(\alpha_0+\beta_0)u_N$. In the case $N=n_0+3$ we
have
\begin{equation*}
\begin{aligned}
u_N&=u_N+(a_1+1)u_{N-1}\\&=\gamma_0u_N+[a_1(\alpha_0+\beta_0)+1]u_{N-1}\\&
=\gamma_0u_N+\gamma_1u_{N-1}.
\end{aligned}
\end{equation*}
In the case $N=n_0+4$ we have
\begin{equation*}
\begin{aligned}
u_N&=u_N+a_1(\gamma_0+1)u_{N-1}+[a_1(a_1+1)+(a_2+1)]u_{N-2}\\
&=\gamma_0u_N+\gamma_1u_{N-1}+[a_1(\alpha_1+\beta_1)+a_2\gamma_0+1]u_{N-2}\\
&=\gamma_0u_N+\gamma_1u_{N-1}+\gamma_2u_{N-2}.
\end{aligned}
\end{equation*}
The next steps follow by induction, and we have recurrence
relation \eqref{AR5} - \eqref{AR6} above.

Therefore,
\begin{equation}\label{AR7}
u_N^2=\left(\sum_{n=n_0+2}^N \gamma_{N-n}u_n\right)^2,
\end{equation}
and similarly to \eqref{ATS7} by Cauchy-Schwartz inequality
\begin{equation}\label{AR8}
\left(\sum_{n=n_0+1}^N \gamma_{N-n}u_n\right)^2\leq
\sum_{n=n_0+2}^N\gamma_{N-n}^2\cdot\sum_{n=n_0+2}^Nu_{n}^{2}
\end{equation}
The equality in \eqref{AR8} is achieved if and only if
$\gamma_{N-n}=cu_n$ for some constant $c$, and since the equality
in \eqref{AR8} is associated with the minimum of the left-hand
side of \eqref{AR8}, the problem reduces to find an appropriate
value $c=c^*$ such that
\begin{equation*}
u_n=c^*\gamma_{N-n}.
\end{equation*}
This finishes the proof.
\end{proof}

The results above are easily extended to general AR($p$) models.
Specifically, we have
\begin{equation}
\label{ARP1}
\widehat{x}_n=a_1\widehat{x}_{n-1}+a_2\widehat{x}_{n-2}+\ldots+a_p\widehat{x}_{n-p}+b,
\end{equation}
and
\begin{equation}
\label{ARP2}
\widetilde{x}_n=a_1\widetilde{x}_{n-1}+a_2\widetilde{x}_{n-2}+\ldots+
a_p\widetilde{x}_{n-p}+b+u_n,
\end{equation}
($n=n_0+p$, $n_0+p+1$+\ldots, $N$), and
\begin{equation}\label{ARP4}
u_N=\widetilde x_{N}-\widehat x_{N}=\sum_{n=n_0+p}^N
\gamma_{N-n}u_n,
\end{equation}
where $\gamma_n=(\alpha_{n,1}+\alpha_{n,2}+\ldots+\alpha_{n,p}),$
and
\begin{eqnarray}\label{ARP5}
\alpha_{0,1} &=& 1,\nonumber\\
\alpha_{0,k} &=& 0, \ k=2,3,\ldots,p,\nonumber\\
\alpha_{n,1} &=& a_1(\alpha_{n-1,1}+\alpha_{n-1,2}+\ldots+\alpha_{n-1,p}) \ \ (n\geq1),\nonumber\\
\alpha_{1,2} &=& 1,\nonumber\\
\alpha_{1,k} &=& 0, \ k=3,4,\ldots,p,\nonumber\\
\alpha_{n,2} &=& a_2(\alpha_{n-2,1}+\alpha_{n-2,2}+\ldots+\alpha_{n-2,p}) \ \ (n\geq2),\nonumber\\
\ldots &=& \ldots\\
\alpha_{k-1,k} &=& 1, \ k=1,2,\ldots,p-1,\nonumber\\
\alpha_{k-1,l} &=& 0, \ l=k+1,k+2,\ldots,p,\nonumber\\
\alpha_{n,k}
&=&
a_k(\alpha_{n-k,1}+\alpha_{n-k,2}\ldots+\alpha_{n-k,p}) \ \ (n\geq k),\nonumber\\
\ldots &=& \ldots\nonumber\\
\alpha_{n,p} &=&
a_p(\alpha_{n-p,1}+\alpha_{n-p,2}\ldots+\alpha_{n-p,p})+1 \ \
(n\geq p).\nonumber
\end{eqnarray}
Thus, similarly to \eqref{AR9} we have the following formula
\begin{equation}\label{ARP9}
\widetilde x_n=a_1\widetilde x_{n-1}+a_2\widetilde
x_{n-2}+\ldots+a_p\widetilde{x}_{n-p}+b+\frac{ \overline{x}-
\widehat x_N}{\sum_{i=n_0+2}^N\gamma_{N-i}^2}\cdot \gamma_{N-n}
\end{equation}
($n=n_0+p, n_0+p+1,\ldots, N-1$), where $\gamma_n$ are now defined
according to \eqref{ARP5}.

\section{Multi-dimensional autoregressive
model}\label{Multi-autoregressive}

In this section we study a multidimensional version of the problem
for AR(1). Let
\begin{equation}\label{MATS1}
\begin{array}{llllllllll}
&{\bf x}_1, &{\bf x}_2, &\ldots, &{\bf x}_{n_0}, &\widetilde{{\bf
x}}_{n_0+1}, &\widetilde{{\bf x}}_{n_0+2}, &\ldots,
&\widetilde{{\bf x}}_{N-1}, &{\bf x}_{N}=\overline{{\bf x}}.
\end{array}
\end{equation}
For this last value we shall also write ${\bf
x}_{N}=\widetilde{\bf x}_{N}$ (with tilde).

As above, the values ${\bf x}_1$, ${\bf x}_2$, \ldots, ${\bf
x}_{n_0}$ are assumed to be observed values, while $\widetilde{\bf
x}_{n_0+1}$, $\widetilde{\bf x}_{n_0+2}$, \ldots, $\widetilde{\bf
x}_{N-1}$ are missing observations.

Taking into consideration only the first $n_0$ observed values one
can build the linear least square predictor as
\begin{equation}
\label{MATS2} \widehat {\bf x}_n=A \widehat {\bf x}_{n-1}+{\bf b}
\end{equation}
for $n=n_0+1, n_0+2,\ldots, N$. Here $A$ is a square matrix, and
${\bf b}$ is a vector.

For $n=n_0+1,\ldots,N$ ($\widetilde {\bf x}_{n_0}$=$\widehat {\bf
x}_{n_0}$) we   find the unknown points by
\begin{equation}\label{MATS3}
\widetilde {\bf x}_n=A \widetilde {\bf x}_{n-1}+{\bf b}+{\bf u}_n.
\end{equation}

The problem is to find the vectors ${\bf u}_n$, $n=n_0+1,\ldots,N$
such that they minimize the sum of squares of their lengths
subject to the constraint that the auto-regression attains the
specified point $\widetilde{\bf x}_{N}$
\begin{equation}
\min_{{\bf u}_n:\; {\bf   x}_{N}=\widetilde{\bf x}_{N}}
\sum_{n=n_0+1}^N \|{\bf u}_n\|^2,\label{MATS4}
\end{equation}
where
\begin{equation*}
\|{\bf u}_n\|=\sqrt{u_{n,1}^2+u_{n,2}^2+\ldots+u_{n,k}^2},
\end{equation*}
and $u_{n,j}$ denotes the $j$th component of the ($k$-dimensional)
vector ${\bf u}_n$.

 According to \eqref{MATS2} and \eqref{MATS3}
\begin{equation}\label{MATS5}
\|{\bf u}_n\|^2=\|\widehat {\bf x}_n-\widetilde {\bf x}_n\|^2,
\end{equation}
and for endpoint ${\bf x}_N$ we have
\begin{equation}\label{MATS6}
\|{\bf u}_N\|^2=\left\|\sum_{n=n_0+1}^NA^{N-n}{\bf
u}_{n}\right\|^2.
\end{equation}
Let $a_{i,j}^{(n)}$ denotes element ($i$, $j$) of matrix $A^n$. We
have the following. The $i$th element of multiplication of
$A^{N-n}$ to vector ${\bf u}_{n}$ can be written as
\begin{equation*}
\sum_{j=1}^ka_{i,j}^{(N-n)}u_{n,j},
\end{equation*}
where $u_{n,j}$ is the $j$th element of the vector ${\bf u}_{n}$.
Therefore \eqref{MATS6} can be written as
\begin{equation}
\label{MATS7}
\begin{aligned}
\|{\bf u}_N\|^2&=\left(\sum_{i=1}^{k}\sum_{n=n_0+1}^{N}\sum_{j=1}^ka_{i,j}^{(N-n)}u_{n,j}\right)^2\\
&=\left(\sum_{n=n_0+1}^{N}\left[\sum_{i=1}^{k}\sum_{j=1}^ka_{i,j}^{(N-n)}u_{n,j}
\right]\right)^2\\
&=\left(\sum_{n=n_0+1}^{N}\sum_{j=1}^{k}\underbrace{\left\{\sum_{i=1}^ka_{i,j}^{(N-n)}\right\}}_{\text{first
term}} \underbrace{u_{n,j}}_{\text{second term}}
\right)^2\\
\end{aligned}
\end{equation}
Therefore, by Cauchy-Schwartz inequality, we have
\begin{equation}
\begin{aligned}
\label{MATS8}\|{\bf
u}_N\|^2&\leq\left(\sum_{n=n_0+1}^{N}\sum_{j=1}^{k}\left[\sum_{i=1}^ka_{i,j}^{(N-n)}\right]^2\right)\times
\left(\sum_{n=n_0+1}^{N}\sum_{j=1}^{k}u_{n,j}^2\right)\\
&=\left(\sum_{n=n_0+1}^{N}\sum_{j=1}^{k}\left[\sum_{i=1}^ka_{i,j}^{(N-n)}\right]^2\right)\times
\left(\sum_{n=n_0+1}^{N}\|{\bf u}_n\|^2\right).
\end{aligned}
\end{equation}
The equality in \eqref{MATS8} is achieved if and only if for some
constant $c$,
\begin{equation}\label{MATS9}
\sum_{i=1}^ka_{i,j}^{(N-n)}=cu_{n,j}
\end{equation}
and similarly to that of Section \ref{autoregressive} the optimal
value of this constant $c^*$ is
\begin{equation}\label{MATS10}
c^*=\frac{\|\widetilde {\bf x}_N- \widehat{\bf
x}_N\|}{\sum_{n=n_0+1}^{N}\sum_{j=1}^{k}\left[\sum_{i=1}^ka_{i,j}^{(N-n)}\right]^2}.
\end{equation}
For the sequence $\mathbf{u}_n$ we have:
\begin{equation}\label{MATS11}
\|{\bf u}_n\|=\frac{\| \overline{\bf x}- \widehat{\bf
x}_N\|}{\sum_{l=n_0+1}^{N}\sum_{j=1}^{k}\left[\sum_{i=1}^k
a_{i,j}^{(N-l)}\right]^2}\cdot \sum_{j=1}^{k}\sum_{i=1}^k
a_{i,j}^{(N-n)},
\end{equation}

Let us now find the vectors ${\bf u}_n$, $n=n_0+1, n_0+2,\ldots,
N$. From \eqref{MATS2} and \eqref{MATS3} we have the following:
\begin{equation}
\label{MATS12} {\bf u}_N=\sum_{n=n_0+1}^{N}A^{N-n}{\bf u}_n.
\end{equation}
Therefore for components of the vector ${\bf u}_N$ we have
equations
\begin{equation}
\label{MATS13} u_{N,i}=\sum_{n=n_0+1}^{N} A_i^{(N-n)}{\bf u}_n,
\end{equation}
where $A_i^{(N-n)}$ denotes the $i$th row of the matrix $A^{N-n}$.
Therefore, by Cauchy-Schwartz inequality
\begin{equation}
\begin{aligned}
\label{MATS14} u_{N,i}^2&=\left(\sum_{n=n_0+1}^{N} A_i^{(N-n)}{\bf
u}_n\right)^2\\
&=\left(\sum_{n=n_0+1}^{N}\sum_{j=1}^{k}a_{i,j}^{(N-n)}
u_{n,j}\right)^2\\
&\leq\left(\sum_{n=n_0+1}^{N}\sum_{j=1}^{k}\left[a_{i,j}^{(N-n)}\right]^2\right)
\left(\sum_{n=n_0+1}^{N}\sum_{j=1}^{k}u_{n,j}^2\right),
\end{aligned}
\end{equation}
where the equality achieves in the case if for some $c_i$
\begin{equation}
\label{MATS15}\sqrt{\sum_{j=1}^{k}\left[a_{i,j}^{(N-n)}\right]^2}=c_i\|{\bf
u}_n\|.
\end{equation}
Therefore, substituting \eqref{MATS15} for \eqref{MATS13} we
obtain:
\begin{equation}
\label{MATS16}c_i=\frac{u_{N,i}}{\sum_{n=n_0}^{N}A_i^{(N-n)}\left[A_i^{(N-n)}\right]^\top}.
\end{equation}

\section{Models of multi-regression}\label{Multi-regression}

Regression models with incomplete data has been studied
intensively in the literature, and there are many approaches the
solution of this problem. The theoretical aspect of the present
approach seems to be new nevertheless.

1. Consider first the following data:
\begin{equation}\label{MMR1}
\begin{array}{lllllllll}
&y_1, &y_2, &\ldots, &y_{n_0}, &\widetilde y_{n_0+1}, &\ldots,
&\widetilde y_{N-1}, &y_N=\overline{y}_N\\
&{\bf x}_1, &{\bf x}_2, &\ldots, &{\bf x}_{n_0}, &{\bf x}_{n_0+1},
&\ldots,
&{\bf x}_{N-1}, &{\bf x}_N.\\
\end{array}
\end{equation}
As above we use the notation $\widetilde y_N=y_N$.

We first find the vector ${\bf a}$ and parameter $b$ by linear
least square predictor, so for $n=n_0+1$, $n_0+2$,\ldots,$N$ we
have
\begin{equation}
\label{MMR2}\widehat{y}_n={\bf a}^\top{\bf x}_n+b.
\end{equation}
We have
\begin{equation*}
\begin{aligned}\widehat y_n&={\bf a}^\top{\bf x}_{n}+b\\
&={\bf a}^\top{\bf x}_{n-1}+b+{\bf a}^\top({\bf x}_n-{\bf x}_{n-1})\\
&=\widehat y_{n-1}+b_n,
\end{aligned}
\end{equation*}
where $b_n={\bf a}^\top({\bf x}_n-{\bf x}_{n-1})$.

Therefore considering
\begin{equation}\label{MMR3}
\widetilde y_n=\widetilde y_{n-1}+b_n+u_n,
\end{equation}
where $\widetilde y_{n_0}=\widehat y_{n_0}$, and $\widetilde
y_N=\overline{y}_N$ and the same problem to
\begin{equation*}
\begin{aligned}
&\text{minimize } \ \sum_{n=n_0+1}^N u_n^2,\\
\end{aligned}
\end{equation*}
we arrive at
\begin{equation}\label{MMR4}
\widetilde y_n=\widetilde y_{n-1}+b_n+\frac{ \overline{y}_N-
\widehat y_N}{N-n_0}.
\end{equation}

2. Let us consider a more extended problem
\begin{equation}\label{MMR5}
\begin{array}{lllllllll}
&{\bf y}_1, &{\bf y}_2, &\ldots, &{\bf y}_{n_0}, &\widetilde {\bf
y}_{n_0+1}, &\ldots,
&\widetilde {\bf y}_{N-1}, &{\bf y}_N=\overline{{\bf y}}_N\\
&{\bf x}_1, &{\bf x}_2, &\ldots, &{\bf x}_{n_0}, &{\bf x}_{n_0+1},
&\ldots,
&{\bf x}_{N-1}, &{\bf x}_N,\\
\end{array}
\end{equation}
where the vectors ${\bf y}$ of the first row all of dimension $m$.

By the linear least square predictor we have
\begin{equation}
\label{MMR6}\widehat{{\bf y}}_n=A{\bf x}_n+{\bf b}.
\end{equation}
Here the vectors $\widehat{{\bf y}}_n$ are of dimension $m$, the
matrix $A$ is of $m\times k$ and the vector ${\bf b}$ is of $m$.
We have:
\begin{equation}
\begin{aligned}
\label{MMR7}\widehat{{\bf y}}_n&=A{\bf x}_n+{\bf b}\\
&=A{\bf x}_{n-1}+{\bf b}+A({\bf x}_n-{\bf x}_{n-1})\\
&=\widehat{{\bf y}}_{n-1}+{\bf b}_n,
\end{aligned}
\end{equation}
where ${\bf b}_n$=$A({\bf x}_n-{\bf x}_{n-1})$.

Let us now consider the equation
\begin{equation}
\label{MMR8}\widetilde{{\bf y}}_n=\widetilde{{\bf y}}_{n-1}+{\bf
b}_n+{\bf u}_n.
\end{equation}
In this specific case we have
\begin{equation}
\label{MMR9}{\bf u}_N=\sum_{n=n_0+1}^{N}{\bf u}_n.
\end{equation}
By the same calculations as earlier (see \eqref{MATS11}) we have:
\begin{equation}
\label{MMR10}\|{\bf u}_n\|=\frac{\|\overline{\bf y}-\widehat{\bf
y}_N\|}{N-n},
\end{equation}
and all the constants $c_i$ defined in Section
\ref{Multi-autoregressive} are the same. Therefore
$u_{n,l}^2=\frac{\|u_n\|^2}{m}$.

We finally have
\begin{equation*}\label{MMR11}
\widetilde {\bf y}_n=\widetilde {\bf y}_{n-1}+{\bf b}_n+\frac{
\overline{\bf y}_N- \widehat {\bf y}_N}{N-n_0}.
\end{equation*}

\section{Numerical work}\label{Numerical work}
Numerical work of this paper consists of two different parts. The
first part is related to the case of interpolating missing data in
autoregressive models. The two numerical results of this part are
reflected in Figure 1a and Figure 1b. The second part of numerical
work is related to two-dimensional autoregressive model. The data
for this model are related to archaeological field and taken from
the paper of Cavanagh, Buck and Litton \cite{Cavanagh Buck and
Litton (1)}.

\subsection{Part 1}
The real data of volatility dynamic of IBM company calculated on
the base of the stock information by the method of
\cite{Goldentayer Klebaner Liptser (1)} have been used for Figure
1a. We removed some data from the middle and the end of this
dynamic and then forecasted missing data by AR(1) model for the
construction of missing data described by \eqref{I2}. The value
$n_0$ is equal to 418, and the corresponding number of missing
data is 85. Then the value $n_1$ is equal to 1058 and the
corresponding number of missing data is 106.

In the second example (Figure 1b) we use the volatility dynamic of
exchange rates of USD and New Israel Shekel. The historical
period, $n_0$, is 1319, and the total length, $N$, is 1466.
Assuming that volatility dynamic is AR(1) model, $v_{n+1}=av_n+b$,
then by calculation of parameters by linear least square predictor
we have $a\approx$0.999576 and $b\approx2.652\cdot 10^{-7}$.
Assuming that volatility dynamics satisfies AR(2) model,
$v_{n+2}=a_1v_{n+1}+a_2v_n+b$, we correspondingly obtain
$a_1\approx2.4$, $a_2\approx-1.4$ and $b\approx2.8\cdot 10^{-7}$.
As we can see, although the difference between these predicted
models is small, both curves are visible in the graph
nevertheless.

\begin{figure}\label{f0}

\includegraphics[width=15cm,height=15cm]{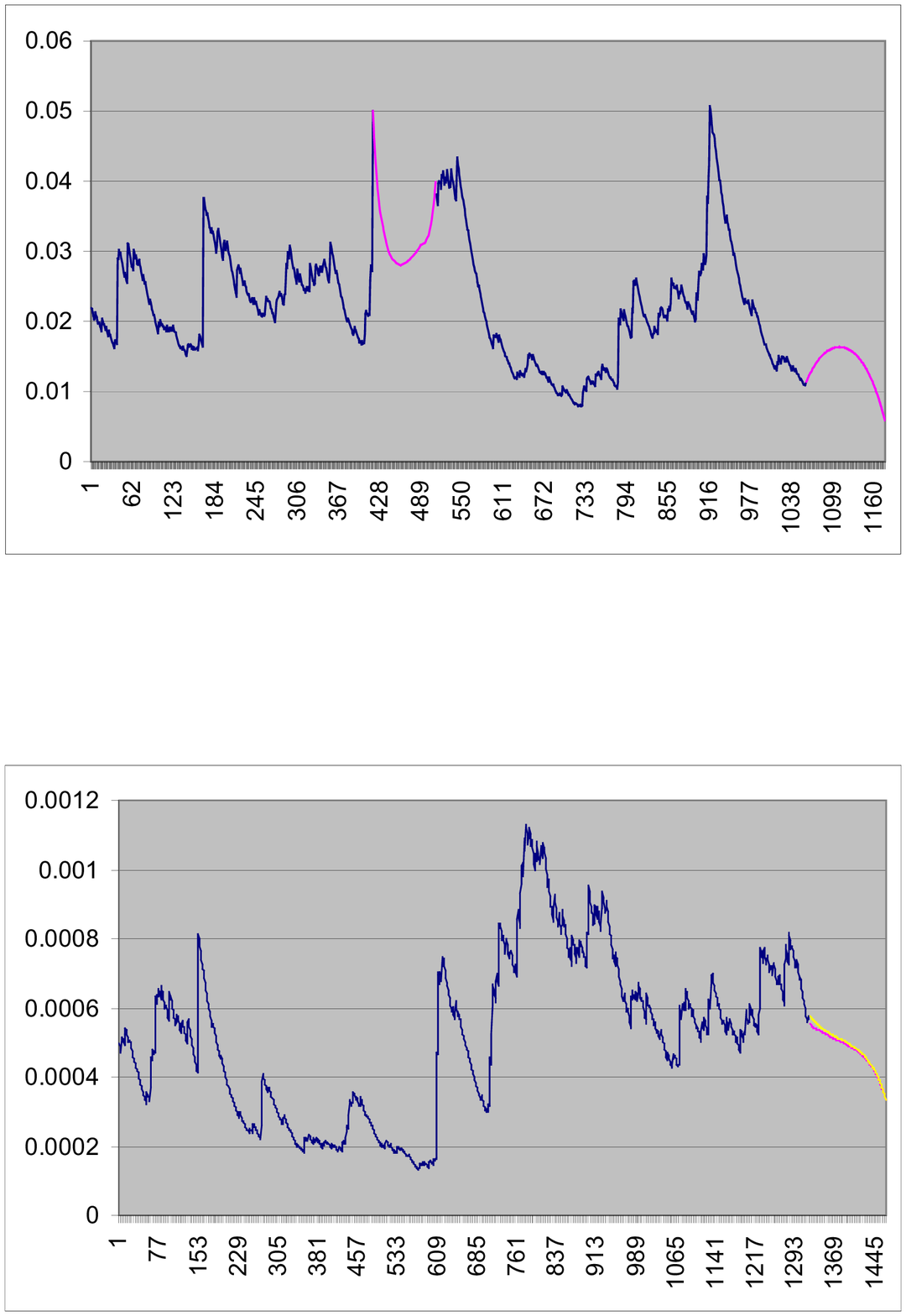}
\caption{
\newline
(a) Forecasting missing data of volatility for IBM Co.: Blue line
- known values, purple line - predicted values by AR(1)
 model;
\newline
 (b) Forecasting missing data of volatility for USD-New
Israel Shekel exchange: Blue line - known values, purple line -
predicted values by AR(1) model and yellow line - predicted values
AR(2) model).}
\end{figure}

\subsection{Part 2} In this part we use data from \cite{Cavanagh Buck and Litton
(1)}. This is data on Phosphate concentration reflected in Figure
1 (p.94). There are missing data in the fifth and sixth row of
these data, and these two rows are the rows of Table 1
corresponding to two-dimensional vector $\mathbf{x}$ with missing
data, where the missing data there are indicated by `+'.
\begin{table}
    \begin{center}
        \begin{tabular}{c|c|c|c|c|c|c|c|c|c|c|c|c|c|c|c}\hline
        59&57&80&71&19&80&60&60&60&62&+&+&166&77&+&68\\
        60&68&75&85&57&44&30&62&38&91&+&+&68&77&+&59\\
        \hline
        \end{tabular}

        \medskip
        \caption{Phosphate concentration (the fragment of data from \cite{Cavanagh Buck and Litton (1)}).}
    \end{center}
\end{table}

We use a first order autoregressive model in order to predict the
missing data. We do not provide all intermediate calculations,
only meaningful results are shown here.

The filling of these missing data is carried out by two steps.
According to our notation $n_0=10$ and $n_1=12$. We have
$\widehat{\mathbf{x}}_{11}=\left(\begin{array}{c}60.43\\
28.57\end{array}\right)$, $\widehat{\mathbf{x}}_{12}=\left(\begin{array}{c}61.10\\
47.98\end{array}\right)$. Next, taking into account the value
$\mathbf{x}_{13}=\left(\begin{array}{c}166\\
68\end{array}\right)$, we obtain the following values for
$\widetilde{\mathbf{x}}_{11}$ and $\widetilde{\mathbf{x}}_{12}$
\begin{equation*}
\widetilde{\mathbf{x}}_{11}=\left(\begin{array}{c}95.10\\
55.33\end{array}\right) \ \ \ \ \ \widetilde{\mathbf{x}}_{12}\left(\begin{array}{c}130.43\\
56.67\end{array}\right).
\end{equation*}
Next, according to the accepted notation, $N$=16. Similarly, we
first
find $\widehat{\mathbf{x}}_{15}$ =$\left(\begin{array}{c}76.90\\
62.93\end{array}\right)$. Then, $\widetilde{\mathbf{x}}_{15}$ =$\left(\begin{array}{c}71.45\\
61.15\end{array}\right)$.

The finally modified table after calculation of missing data is
now Table 2.
\begin{table}
    \begin{center}
        \begin{tabular}{c|c|c|c|c|c|c|c|c|c|c|c|c|c|c|c}\hline
        59&57&80&71&19&80&60&60&60&62&95.10&130.43&166&77&71.45&68\\
        60&68&75&85&57&44&30&62&38&91&55.33&56.67&68&77&61.15&59\\
        \hline
        \end{tabular}

        \medskip
        \caption{Phosphate concentration (the finally modified table).}
    \end{center}
\end{table}

\end{document}